\documentclass[11pt]{article}
\usepackage{amsmath,amssymb,eucal}
\usepackage[dvips]{graphicx}

\newtheorem{thm}{\bf Theorem}[section]
\newtheorem{co}[thm]{\bf Corollary}
\newtheorem{lm}[thm]{\bf Lemma}
\newtheorem{rem}[thm]{\bf Remark}
\newtheorem{prop}[thm]{\bf Proposition}
\newtheorem{ex}[thm]{\bf Example}

\newtheorem{qu}[thm]{\bf Question}

\numberwithin{equation}{section}

\def\qed{$\Box$}

\begin{document}

\title{{\bf 
{\large Hyperbolicity and identification of Berge knots of types VII and VIII}}}
\author{{\normalsize Teruhisa KADOKAMI}}
\date{{\normalsize December 21, 2011}}
\footnotetext[0]{%
}
\maketitle

\begin{abstract}{
T.~Saito and M.~Teragaito asked whether Berge knots of type VII 
are hyperbolic,
and showed that some infinite sequences of the knots are hyperbolic.
We show that Berge knots of types VII and VIII are hyperbolic
except the known sequence of torus knots.
We used the Reidemeister torsions.
As a result, the Alexander polynomials of them
have already shown their hyperbolicities.
We also show that the standard parameters
identify Berge knots of types VII and VIII, and study
what kind of information identify them.
}\end{abstract}

\section{Introduction}~\label{sec:intro}
\begin{figure}[ht]
\begin{center}
\includegraphics[scale=0.56]{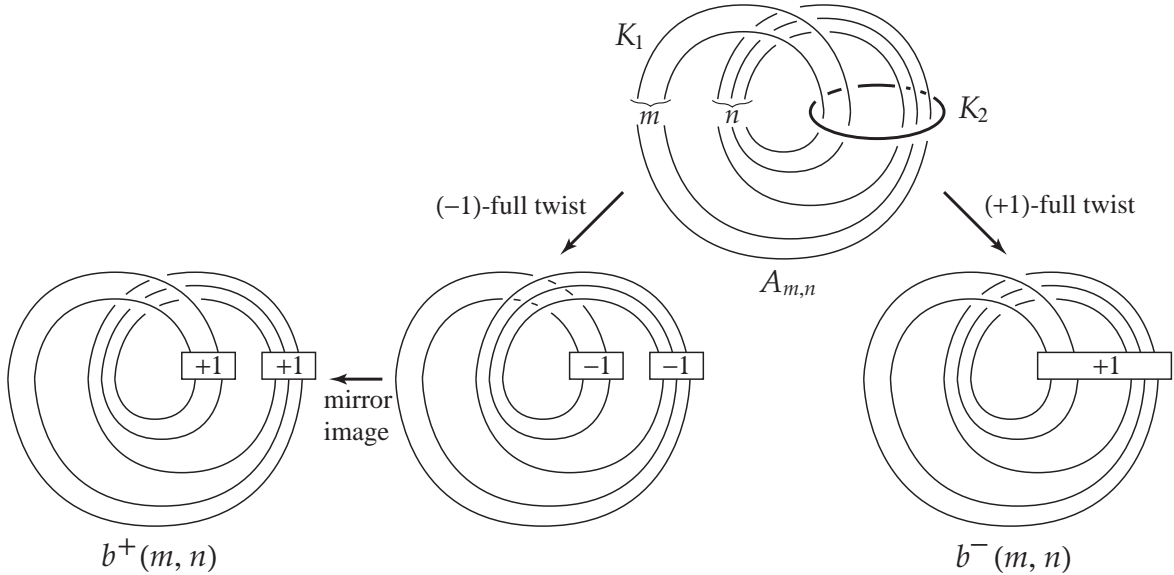}
\caption{$A_{m,n}$, and Berge knots of type VII ($b^+(m, n)$) 
and VIII ($b^-(m, n)$)}
\label{fig:k+-}
\end{center}
\end{figure}

For two positive coprime integers $m$ and $n$,
let $A_{m,n}=K_1\cup K_2$ be a $2$-component link 
in the upper side of Figure 1, where
$K_1$ is an $(m, n)$-torus knot $T(m, n)$
and $K_2$ is the unknot.
By operating $1$-surgery (i.e.\ $(-1)$-full twist) along $K_2$
and taking the mirror image, we obtain a knot $b^+(m, n)$
in the lefthand side of Figure 1.
By operating $(-1)$-surgery (i.e.\ $(+1)$-full twist) along $K_2$, 
we obtain a knot $b^-(m, n)$ in the righthand side of Figure 1.
Then the set of Berge knots of types VII  \cite{Ber}
(the set of Berge knots of types VIII, respectively)
is the set of $\{b^+(m, n)\}$ and their mirror images
(the set of $\{b^-(m, n)\}$ and their mirror images, respectively).
Berge knots of types VII  and VIII are characterized
that they can be situated on the fiber surfaces of
the trefoil and the figure eight knot respectively.
We say that the {\it standard parameters} of 
$b^+(m, n)$ and $b^-(m, n)$ are
$(+, m, n)$ and $(-, m, n)$, respectively.
The notations in \cite{Yam1} (see also \cite{KY, ST, Yam2, Yam3})
are $k^+(m, n)=b^+(m, n)$ and $k^-(m, m+n)=b^-(m, n)$,
respectively.
Since $A_{m,n}=A_{n,m}$ as ordered oriented links,
we have $b^+(m, n)=b^+(n, m)$ and $b^-(m, n)=b^-(n, m)$
as oriented knots, and we may assume $1\le m\le n$ and $\gcd(m, n)=1$.
As trivial cases, $b^+(1, n)$ and $b^-(1, n)$ are 
an $(n+1, n)$-torus knot and an $(n+1, n+2)$-torus knot respectively,
and we have $b^+(1, n)=b^-(1, n-1)=T(n, n+1)$.
Moreover since $A_{m,n}$ is invertible,
both $b^+(m, n)$ and $b^-(m, n)$ are also invertible.
Hence we study hyperbolicity and identification problem of 
$b^+(m, n)$ and $b^-(m, n)$ for $2\le m<n$.

\medskip

T.~Saito and M.~Teragaito \cite{ST} asked 
hyperbolicity of $b^+(m, n)$,
and showed that some infinite sequences of the knots
are hyperbolic by using 
the necessary and sufficient condition for hyperbolicity 
of a Berge knot due to Saito \cite{Sai}, and
by elementary number theoretical results on Fibonacci series.
In the present paper, 
we show that both $b^+(m, n)$ and $b^-(m, n)$
are hyperbolic for $2\le m<n$ (Theorem \ref{thm:hyp}).
We used the Reidemeister torsions.
From the proof of Theorem \ref{thm:hyp}
and the author's previous works \cite{Kad1, Kad2},
the Alexander polynomials of them have already shown their hyperbolicities 
(Theorem \ref{thm:Alex}).

\medskip

In \cite[Theorem 1.1]{ST}, it is shown that there are infinite pairs 
of two distinct knots in $S^3$ whose same integer surgeries
yield homeomorphic lens spaces (Theorem \ref{thm:ST}).
We point out that the resulting lens space up to orientations by 
an odd integer surgery along a hyperbolic Berge knot of type VII or VIII
identifies the standard parameter and the knot itself 
up to mirror images 
(Lemma \ref{lm:lens}, Corollary \ref{co:mirror}).
We show that a pair of the odd surgery coefficient and 
the degree of the Alexander polynomial,
which is equivalent to the Seifert genus in this case,
identifies a Berge knot of type VII or VIII up to mirror images
(Theorem \ref{thm:detect}, Corollary \ref{co:pg}).
Y.~Yamada \cite{Yam3} defined $b^+(m, n)$ and $b^-(m, n)$
for every coprime pair $(m, n)$, and discussed equivalence 
among them up to orientations and mirror images.
For every element of the extended class, 
there is one representative with $1\le m\le n$ 
(see Lemma \ref{lm:primitive} (3)).
Since a double torus knot is strongly invertible, 
a non-trivial torus knot is non-amphicheiral, and
a non-torus knot yielding a lens space is non-amphicheiral
by the Cyclic Surgery Theorem \cite{CGLS}, 
the extended class can be also classified completely
(see also Remark \ref{rem:Yam} (2)).

\medskip

We ask whether one of
the odd surgery coefficient $p$ and 
the half of the degree of the Alexander polynomial $g$
of a Berge knot of type VII or VIII determines uniquely
the standard parameter or not.
For the case of $p$, we give the complete answer 
which has been already described in \cite{Ber} 
by a theory of quadratic fields (Lemma \ref{lm:prime}).
We make a complete table of $p$, $(\pm, m, n)$ and $g$ for $p\le 500$.
By the table, we can give examples for the case of $g$.

\medskip

K.~Baker \cite{Ba} showed that
families of Berge knots of types VII and VIII include 
hyperbolic knots with arbitrarily large volume.
Recently J.~Greene \cite{Gr} showed that
every doubly primitive knot is a Berge knot.

\medskip

In Section \ref{sec:surgery},
we give a definition of the Reidemeister torsion, and
prepare the surgery formulae for the Reidemeister torsions
and the Franz lemma.
In Section \ref{sec:Alex},
we compute the Alexander polynomials of 
$b^+(m, n)$ and $b^-(m, n)$ from that of $A_{m,n}$,
determine their degrees, and 
compute the Reidemeister torsions of lens spaces
obtained from $b^+(m, n)$ and $b^-(m, n)$.
In Section \ref{sec:hyp},
we show that both $b^+(m, n)$ and $b^-(m, n)$
are hyperbolic for $2\le m<n$ (Theorem \ref{thm:hyp}).
In Section \ref{sec:detect},
we give the complete answer for identification problem
of $b^+(m, n)$ and $b^-(m, n)$ (Theorem \ref{thm:detect}),
and make a table of the standard parameters.
In Section \ref{sec:rem}, we give a remark about
a relation between our method and 
K.~Ichihara, T.~Saito and M.~Teragaito \cite{IST}.

\medskip

Throughout the paper, for coprime integers $p$ and $q$,
a lens space of type $(p, q)$ is the result of $-p/q$-surgery 
along the unknot, and it is denoted by $L(p, q)$.

\section{Surgery formulae for the Reidemeister torsions, and Franz lemma}
\label{sec:surgery}
For a precise definition of the Reidemeister torsion,
the reader refer to V.~Turaev \cite{Tur}.
In this section, we give surgery formulae of the Reidemeister torsions,
and the Franz lemma.
Throughout the paper, $\zeta_d$ denotes a primitive $d$-th root of unity,
and $(\mathbb{Z}/d\mathbb{Z})^{\times}$ is the multiplicative group
of a ring $\mathbb{Z}/d\mathbb{Z}$.

\medskip

Let $X$ be a finite CW complex
and $\pi : {\tilde X}\to X$ its maximal abelian covering.
Then ${\tilde X}$ has a CW structure induced by that of $X$ and $\pi$,
and the cell chain complex $\mathbf{C}_{\ast}$ of ${\tilde X}$
has a $\mathbb{Z}[H]$-module structure, where
$H=H_1(X ; \mathbb{Z})$ is the first homology of $X$.
For an integral domain $R$ and a ring homomorphism
$\psi : \mathbb{Z}[H]\to R$,  
$\mathbf{C}_{\ast}^{\psi}
=\mathbf{C}_{\ast}\otimes_{\mathbb{Z}[H]}Q(R)$
is the chain complex of ${\tilde X}$ related with $\psi$,
where $Q(R)$ is the quotient field of $R$.
Note that there is a natural inclusion $R\hookrightarrow Q(R)$.
Then the Reidemeister torsion of $X$ related with $\psi$, 
denoted by $\tau^{\psi}(X)$,
is obtained from $\mathbf{C}_{\ast}^{\psi}$, and 
is an element of $Q(R)$ up to multiplication of $\pm \psi(h)\ (h\in H)$.
If $R=\mathbb{Z}[H]$ and $\psi$ is the identity map,
then we denote $\tau^{\psi}(X)$ by $\tau(X)$.
We note that $\tau^{\psi}(X)$ is not zero if and only if 
$\mathbf{C}_{\ast}^{\psi}$ is acyclic.
If $A$ and $B$ are elements of $Q(R)$, and there exists an element $h\in H$
such that $A=\pm \psi(h)B$, then we denote the equality by $A\doteq B$.

\medskip

Let $H$ be a finitely generated abelian group
with a direct sum $H\cong B\oplus T$,
where $B$ is a free abelian group of rank $r$
and $T$ is a finite abelian group.
Then, by the Chinese Remainder Theorem or the Pontrjagin duality,
$\mathbb{Q}[H]$ is a direct sum, as commutative rings,
of $\mathbb{Q}(\zeta_{d_i})\otimes F\cong F(\zeta_{d_i})$,
where $F=\mathbb{Q}(B)$ is a rational function field with $r$ variables 
over $\mathbb{Q}$.
In particular, 
if $H=\langle t\ |\ t^p=1\rangle \cong \mathbb{Z}/p\mathbb{Z}$, then
we have the canonical isomorphism:
\begin{equation}\label{eq:sum1}
\mathbb{Q}[t, t^{-1}]/(t^p-1)\cong 
\bigoplus_{d|p, d\ge 1}\mathbb{Q}[t, t^{-1}]/(\mathbf{\Phi}_d(t)),
\end{equation}
where $\mathbf{\Phi}_d(t)$ is the $d$-th cyclotomic polynomial, and
$t\ (\mathrm{mod}\ \! t^p-1)$ in the lefthand side corresponds to
$(t\ (\mathrm{mod}\ \! \mathbf{\Phi}_d(t)))_{d|p, d\ge 1}$
in the righthand side.
(\ref{eq:sum1}) deduces an isomorphism:
\begin{equation}\label{eq:sum2}
\mathbb{Q}[H]\cong \bigoplus_{d|p, d\ge 1}\mathbb{Q}(\zeta_d).
\end{equation}
However isomorphisms
\begin{equation}\label{eq:noncano}
\mathbb{Q}[t, t^{-1}]/(t^p-1)\cong \mathbb{Q}[H]
\quad \mbox{and}\quad
\mathbb{Q}[t, t^{-1}]/(\mathbf{\Phi}_d(t))\cong \mathbb{Q}(\zeta_d)
\end{equation}
are not uniquely determined.
In the former case, $t^i$ can be a generator of $H$
for any $i\in (\mathbb{Z}/p\mathbb{Z})^{\times}$.
In the latter case, $t$ can correspond to $\zeta_d^i$
for any $i\in (\mathbb{Z}/d\mathbb{Z})^{\times}$.
In discussing the Reidemeister torsions,
only the field components in the direct sum of $\mathbb{Q}[H]$
are essential.

\medskip

Let $\psi_d : \mathbb{Z}[t, t^{-1}]/(t^p-1)\to \mathbb{Q}(\zeta_d)$
be a ring homomorphism defined by
$\psi_d(t)=\zeta_d$.
The map $\psi_d$ can be regarded as a projection to 
a field component in (\ref{eq:sum2}).
Since a primitive $d$-th root of unity is not unique for $d\ge 3$,
the map $\psi_d$ has some possibilities.
Its ambiguity is the Galois group of $\mathbb{Q}(\zeta_d)$
over $\mathbb{Q}$, which is isomorphic to
$(\mathbb{Z}/d\mathbb{Z})^{\times}$.

\medskip

Let $E$ be a compact $3$-manifold whose boundary $\partial E$ 
consists of tori, $M=E\cup V$ a $3$-manifold 
obtained by attaching a solid torus $V$ to $E$ by an attaching map
$f : \partial V\to \partial E$, and
$\iota : E\hookrightarrow M$ the natural inclusion.
Let $l'$ be the core of $V$, and
$[l']$ its homology class in $H_1(M)$.

\begin{lm}\label{lm:surgery1}
{\rm (Surgery formula I)}
Let $F$ be a field and $\psi : \mathbb{Z}[H_1(M)]\to F$ a ring homomorphism.
If $\psi([l'])\ne 1$, then
$$\tau^{\psi}(M)\doteq 
\tau^{\psi'}(E)\cdot (\psi([l'])-1)^{-1},$$
where $\psi'=\psi \circ \iota_{\ast}$ 
($\iota_{\ast}$ is a ring homomorphism induced by $\iota$).
\end{lm}

For an oriented $\mu$-component link 
$L=K_1\cup \ldots \cup K_{\mu}$ in $S^3$, 
let $E_L$ be the complement of $L$, and
${\it \Delta}_L(t_1, \ldots, t_{\mu})$ the Alexander polynomial of $L$
where $t_i$ is represented by a meridian of $K_i$.
The Reidemeister torsion is closely related with the Alexander polynomial.

\begin{lm}\label{lm:Alexander}
{\rm (Milnor \cite{Mil})}
We have
$$\tau(E_L)\doteq
\left\{
\begin{array}{cl}
{\it \Delta}_L(t_1)(t_1-1)^{-1} & (\mu =1),\\
{\it \Delta}_L(t_1, \ldots, t_{\mu}) & (\mu \ge 2).
\end{array}
\right.$$
\end{lm}

For an oriented $\mu$-component link 
$L=K_1\cup \ldots \cup K_{\mu}$ in $S^3$, 
let $m_i$ and $l_i$ be a meridian and a longitude 
of the $i$-th component $K_i$,
$[m_i]$ and $[l_i]$ their homology classes, and
$(L ; p_1/q_1, \ldots, p_{\mu}/q_{\mu})$ 
the result of $(p_1/q_1, \ldots, p_{\mu}/q_{\mu})$-surgery along $L$
where $p_i$ and $q_i$ are coprime integers or $p_i/q_i=\emptyset$.
We take integers $r_i$ and $s_i$ satisfying $p_i s_i -q_i r_i = -1$.
For coprime integers $p$ and $q$, 
the inverse of $q$ in $(\mathbb{Z}/p\mathbb{Z})^{\times}$
is denoted by $\overline{q} \ (\mathrm{mod}\ \! p)$, and also 
$\overline{q}\ (\mathrm{mod}\ \! d)$ for any divisor $d$ ($\ge 2$) of $p$.

\medskip

\begin{lm}\label{lm:surgery2}
{\rm (Surgery formula II;\ Sakai \cite{Sak}, Turaev \cite{Tur})}
\begin{enumerate}
\item[(1)]
For $M=(K ; p/q)\ (|p|\ge 2)$
and a diviser $d\ge 2$ of $p$, 
we have 
$$\tau^{\psi_d}(M)\doteq {\mit \Delta}_K(\zeta_d)
(\zeta_d-1)^{-1}(\zeta_d^{{\bar q}}-1)^{-1}$$
where $q\overline{q}\equiv 1\ (\mathrm{mod}\ \! p)$.

\item[(2)]
For $M=(L ; p_1/q_1, \ldots, p_{\mu}/q_{\mu})\ (\mu \ge 2)$,
let $F$ be a field and $\psi : \mathbb{Z}[H_1(M)]\to F$ a ring homomorphism.
If $\psi([m_i]^{r_i}[l_i]^{s_i})\ne 1\ (i=1, \ldots, \mu)$, then we have
$$\tau^{\psi}(M)\doteq {\it \Delta}_L(\psi([m_1]), \ldots, \psi([m_{\mu}]))
\prod_{i=1}^{\mu}(\psi([m_i]^{r_i}[l_i]^{s_i})-1)^{-1}.$$
\end{enumerate}
\end{lm}

\begin{ex}\label{ex:lens}
{\rm
By Lemma \ref{lm:surgery2} (1), for a divisor $d\ge 2$ of $p$, 
we have
$$\tau^{\psi_d}(L(p, q))\doteq (\zeta_d-1)^{-1}(\zeta_d^{{\bar q}}-1)^{-1}$$
where $q\overline{q}\equiv 1\ (\mathrm{mod}\ \! p)$.}
\end{ex}

\begin{rem}\label{rem:amb}
{\rm
By the ambiguity in (\ref{eq:noncano}), 
the Reidemeister torsion for the result of Dehn surgery along a knot
can include information of a meridian of the knot.}
\end{rem}

W.~Franz \cite{Fz} showed the following, and
classified lens spaces by using it.

\begin{lm}\label{lm:Franz}
{\rm (Franz \cite{Fz})}
For $a_i$ and $b_i\in (\mathbb{Z}/p\mathbb{Z})^{\times}$\ $(i=1, \ldots, n)$,
suppose
$$
\prod_{i=1}^n(\zeta_p^{a_i}-1)
\doteq \prod_{i=1}^n(\zeta_p^{b_i}-1).
$$
Then there exists a permutation $\sigma$ of $\{1, \ldots, n\}$
such that $a_i=\pm b_{\sigma(i)}$ for all $i=1, \ldots, n$.
In other words, 
$\{ \pm a_i \, (\mathrm{mod}\ \! p)\} = \{ \pm b_i \, \, (\mathrm{mod}\ \! p)\}$
as multiple sets.
\end{lm}

\section{The Alexander polynomials of $b^+(m, n)$ and $b^-(m, n)$}
\label{sec:Alex}
We compute the Alexander polynomials of 
$b^+(m, n)$ and $b^-(m, n)$ from that of $A_{m,n}$,
determine their degrees, and 
compute the Reidemeister torsions of lens spaces
obtained from $b^+(m, n)$ and $b^-(m, n)$.
Some results are from 
a work of the author with Y.~Yamada \cite{KY}.

\medskip

For a coprime positive pair $(m, n)$, 
we define a set ${\frak I}(m, n)$ by 
\begin{eqnarray*}
{\frak I}(m, n) 
& = & 
(m\mathbb{Z} \cup n\mathbb{Z}) \cap 
\{ k \in \mathbb{Z}\ \vert \ 0 \leq k \leq mn \} \\
& = &
\{0, m, 2m, \ldots, nm \} \cup 
\{0, n, 2n, \ldots, mn \}
\end{eqnarray*}
Note that the cardinality of ${\frak I}(m, n)$ is $m+n$.
We sort the all elements in ${\frak I}(m, n)$ as 
$$0 = k_0 < k_1 < k_2 < \cdots < k_{m+n-1} = mn \qquad
(k_i \in {\frak I}(m, n)).$$

\begin{thm}\label{thm:AmnAlex}
{\rm (\cite{KY})}
The Alexander polynomial of the link $A_{m,n}$ is 
$${\it \Delta}_{A_{m,n}} (t,x)\doteq
\sum_{i=0}^{m+n-1} t^{k_i} x^i.$$
\end{thm}

We define $u_j\ (j=0, 1, 2, \ldots, n)$ and
$w_j\ (j=0, 1, 2, \ldots, m)$ by
$$k_{u_j}=jm\quad \mbox{and}\quad k_{w_j}=jn.$$
Then we have the following:

\begin{lm}\label{lm:Imn}
{\rm (\cite{KY})}
\begin{enumerate}
\item[(1)]
$$u_j=\left\{
\begin{array}{cl}
0 & (j=0),\medskip\\
\left[\frac{jm}{n}\right]+j &
(j=1, 2, \ldots, n-1),\medskip\\
m+n-1 & (j=n),
\end{array}
\right.$$
and $u_j+u_{n-j}=m+n-1$,
where $[\ \cdot\ ]$ is the gaussian symbol.

\item[(2)]
$$w_j=\left\{
\begin{array}{cl}
0 & (j=0),\medskip\\
\left[\frac{jn}{m}\right]+j &
(j=1, 2, \ldots, m-1),\medskip\\
m+n-1 & (j=m),
\end{array}
\right.$$
and $w_j+w_{m-j}=m+n-1$.
\end{enumerate}
\end{lm}

By Lemma \ref{lm:Imn} (1) and (2), we have
$$\sum_{i=0}^{m+n-1}t^{k_i}x^i
=1+\sum_{i=1}^{n-1}t^{im}x^{u_i}+\sum_{j=1}^{m-1}t^{jn}x^{w_j}
+t^{mn}x^{m+n-1}.$$

Here we denote $b^+(m, n)$ and $b^-(m, n)$ by $K^+$ and $K^-$, respectively. 
Let $g$ be the half of the degree of the Alexander polynomial of $K^{\pm}$.
Since every knot yielding a lens space is fibered Y.~Ni \cite{Ni}, 
the Seifert genus of $K^{\pm}$ is $g$.

\begin{prop}\label{prop:genus}
{\rm (Saito and Teragaito \cite[Lemma 5.1]{ST})}
The Alexander polynomial of $K^{\pm}$ is
$${\it \Delta}_{K^{\pm}}(t)\doteq
\frac{t-1}{t^{m+n}-1}\cdot
\sum_{i=0}^{m+n-1}t^{k_i\mp i(m+n)},$$
and the half of the degree of the Alexander polynomial of $K^{\pm}$ is
$$g=\frac{(m+n-1)^2\mp mn}{2}.$$
\end{prop}

\noindent
{\bf Proof}\ 
For $A_{m,n}=K_1\cup K_2$ and
$M=(A_{m,n}; \emptyset, \pm 1)$, 
let $E$ be the complement of $A_{m,n}$,
$V$ an attaching solid torus to $\partial N(K_2)\subset \partial E$ 
where $N(K_2)$ is a regular neighborhood of $K_2$ in $S^3$,
and $M=E\cup V$.
Let $m_i$ and $l_i$ $(i=1, 2)$ be a meridian and a longitude of $K_i$ respectively,
and $m'$ and $l'$ be a meridian and a core of $V$ respectively.
We note that $M$ is homeomorphic to the complement of $K^{\pm}$.
Then we study homological conditions.

\medskip

In $H_1(E)$, 
\begin{equation*}
[l_1]=[m_2]^{m+n},\quad [l_2]=[m_1]^{m+n},
\end{equation*}
and
\begin{eqnarray*}
H_1(E) & \cong & 
\langle [m_1], [l_1], [m_2], [l_2]\ |\ 
[l_1]=[m_2]^{m+n}, [l_2]=[m_1]^{m+n}\rangle \\
& \cong & 
\langle [m_1], [m_2]\ |\ - \rangle \cong \mathbb{Z}^2.
\end{eqnarray*}
In $H_1(M)$, 
\begin{equation*}
[m']=[m_2]^{\pm 1}[l_2]=1,\quad [l']=[m_2],
\end{equation*}
and
\begin{eqnarray*}
H_1(M) & \cong & 
\langle [m_1], [m_2]\ |\ 
[m_1]^{m+n}[m_2]^{\pm 1}=1\rangle \\
& \cong & 
\langle t\ |\ - \rangle \cong \mathbb{Z},
\end{eqnarray*}
where $t=[m_1]$.
Then
\begin{equation}\label{eq:rel0}
\begin{matrix}
[m_2] & = & t^{\mp (m+n)}, \hfill \medskip \\
[l'] & = & [m_2]=t^{\mp (m+n)}. \hfill
\end{matrix}
\end{equation}

By Theorem \ref{thm:AmnAlex}, (\ref{eq:rel0}) and
the surgery formula I (Lemma \ref{lm:surgery1}),
we have the Alexander polynomial of $K^{\pm}$ as
\begin{equation*}
{\displaystyle
{\it \Delta}_{K^{\pm}}(t)\doteq
\frac{(t-1){\it \Delta}_{A_{m,n}}(t, t^{\mp (m+n)})}
{t^{m+n}-1}\doteq 
\frac{t-1}{t^{m+n}-1}\cdot
\sum_{i=0}^{m+n-1}t^{k_i\mp i(m+n)}}
\end{equation*}

\noindent
(i)\ The $K^+$ case

\medskip

If $k_i=jm\ (j=0, 1, 2, \ldots, n-1)$ (i.e.\ $i=u_j$), then
by Lemma \ref{lm:Imn} (1), we have
$$k_i-i(m+n)
=jm-\left( j+\left[ \frac{jm}{n}\right]\right)(m+n)
=-\left\{ jn+\left[ \frac{jm}{n}\right](m+n)\right\},$$
and it is monotonely decreasing about $j$.
Moreover
\begin{eqnarray*}
& & k_{u_{n-1}}-(m+n-2)(m+n)=(n-1)m-(m+n-2)(m+n)
\medskip\\
& > &
k_{u_n}-(m+n-1)(m+n)=-\{(m+n-1)(m+n)-mn\}.
\end{eqnarray*}

If $k_i=jn\ (j=0, 1, 2, \ldots, m-1)$ (i.e.\ $i=w_j$), then
by Lemma \ref{lm:Imn} (2), we have
$$k_i-i(m+n)
=jn-\left( j+\left[ \frac{jn}{m}\right]\right)(m+n)
=-\left\{ jm+\left[ \frac{jn}{m}\right](m+n)\right\},$$
and it is monotonely decreasing about $j$.
Moreover
\begin{eqnarray*}
& & k_{w_{m-1}}-\left(m+n-2-\left[ \frac nm\right]\right)
(m+n)=(m-1)n-\left(m+n-2-\left[ \frac nm\right]\right)(m+n)
\medskip\\
& > &
k_{w_m}-(m+n-1)(m+n)=-\{(m+n-1)(m+n)-mn\}.
\end{eqnarray*}

Hence the degree of ${\it \Delta}_{K^+}(t)$ is
$$2g=\{(m+n-1)(m+n)-mn\}+1-(m+n)
=(m+n-1)^2-mn.$$

\noindent
(ii)\ The $K^-$ case

\medskip

If $k_i=jm\ (j=0, 1, 2, \ldots, n-1)$ (i.e.\ $i=u_j$), then
by Lemma \ref{lm:Imn} (1), we have
$$k_i+i(m+n)
=jm+\left( j+\left[ \frac{jm}{n}\right]\right)(m+n)$$
is monotonely increasing about $j$.
Moreover
\begin{eqnarray*}
& & k_{u_{n-1}}+(m+n-2)(m+n)=(n-1)m+(m+n-2)(m+n)
\medskip\\
& < &
k_{u_n}+(m+n-1)(m+n)=mn+(m+n-1)(m+n).
\end{eqnarray*}

If $k_i=jn\ (j=0, 1, 2, \ldots, m-1)$ (i.e.\ $i=w_j$), then
by Lemma \ref{lm:Imn} (2), we have
$$k_i+i(m+n)
=jn+\left( j+\left[ \frac{jn}{m}\right]\right)(m+n)$$
is monotonely increasing about $j$.
Moreover
\begin{eqnarray*}
& & k_{w_{m-1}}+\left(m+n-2-\left[ \frac nm\right]\right)
(m+n)=(m-1)n+\left(m+n-2-\left[ \frac nm\right]\right)(m+n)
\medskip\\
& < &
k_{w_m}+(m+n-1)(m+n)=mn+(m+n-1)(m+n).
\end{eqnarray*}

Hence the degree of ${\it \Delta}_{K^-}(t)$ is
$$2g=\{(m+n-1)(m+n)+mn\}+1-(m+n)
=(m+n-1)^2+mn.$$
\qed

\bigskip

The following lemma is obtained from
Kirby moves \cite{Kir} and Rolfsen moves \cite{Ro}.

\begin{lm}\label{lm:Amnlens}
{\rm (\cite{KY, Yam3})}
$M=(A_{m,n};mn, r)$ is a lens space for every integer $r\in \mathbb{Z}$.
Then $M\cong L(rmn-(m+n)^2, m\overline{n})$.
In particular, we have
$$\mp (A_{m,n};mn, \pm 1)\cong (K^{\pm}; (m+n)^2\mp mn)
\cong L((m+n)^2\mp mn, m\overline{n}).$$
\end{lm}

For $A_{m,n}=K_1\cup K_2$ and
$M=(A_{m,n};mn, r)$ with $r\in \mathbb{Z}$, 
let $E$ be the complement of $A_{m,n}$,
$V_i$ $(i=1, 2)$ an attaching solid torus to $\partial N(K_i)\subset \partial E$ 
where $N(K_i)$ is a regular neighborhood of $K_i$ in $S^3$,
and $M_1=E\cup V_1$ (then $M=M_1\cup V_2$).
Let $m_i$ and $l_i$ $(i=1, 2)$ be a meridian and a longitude of $K_i$ respectively,
and $m_i'$ and $l_i'$ $(i=1, 2)$ be a meridian and a core of $V_i$ respectively.
We set $p= (m+n)^2-rmn$.
Then we study homological conditions.

\medskip

In $H_1(E)$, 
\begin{equation*}
[l_1]=[m_2]^{m+n},\quad [l_2]=[m_1]^{m+n},
\end{equation*}
and
\begin{eqnarray*}
H_1(E) & \cong & 
\langle [m_1], [l_1], [m_2], [l_2]\ |\ 
[l_1]=[m_2]^{m+n}, [l_2]=[m_1]^{m+n}\rangle \\
& \cong & 
\langle [m_1], [m_2]\ |\ - \rangle \cong \mathbb{Z}^2.
\end{eqnarray*}
In $H_1(M_1)$, 
\begin{equation*}
[m_1']=[m_1]^{mn}[l_1]=1,\quad [l_1']=[m_1],
\end{equation*}
and
\begin{eqnarray*}
H_1(M_1) & \cong & 
\langle [m_1], [m_2]\ |\ 
[m_1]^{mn}[m_2]^{m+n}=1\rangle \\
& \cong & 
\langle T\ |\ - \rangle \cong \mathbb{Z},
\end{eqnarray*}
where $T=[m_1]^u[m_2]^v$ for integers $u$ and $v$
such that $(m+n)u-mnv=1$.
Then
\begin{equation}\label{eq:rel1}
\begin{matrix}
[m_1] & = & [m_1]^{(m+n)u-mnv}
=([m_1]^u[m_2]^v)^{m+n}([m_1]^{mn}[m_2]^{m+n})^{-v}
=T^{m+n}, \hfill \medskip \\
[m_2] & = & [m_2]^{(m+n)u-mnv}
=([m_1]^u[m_2]^v)^{-mn}([m_1]^{mn}[m_2]^{m+n})^u
=T^{-mn}, \hfill \medskip \\
[l_1'] & = & [m_1]=T^{m+n}. \hfill
\end{matrix}
\end{equation}
In $H_1(M)$, 
\begin{equation*}
[m_2']=[m_2]^r[l_2]=1,\quad [l_2']=[m_2],
\end{equation*}
and
\begin{eqnarray*}
H_1(M) & \cong & 
\langle [m_1], [m_2]\ |\ 
[m_1]^{mn}[m_2]^{m+n}=1, 
[m_1]^{m+n}[m_2]^r=1\rangle \\
& \cong & 
\langle T\ |\ T^p=1 \rangle \cong \mathbb{Z}/p\mathbb{Z},
\end{eqnarray*}
where $T$ is the induced one from $T$ in (\ref{eq:rel1}).
Then
\begin{equation}\label{eq:rel2}
[l_2']=[m_2]=T^{mn}.
\end{equation}

\begin{lm}\label{lm:formula}
{\rm (\cite{KY})}
Let $T$ be an indeterminancy.
$${\it \Delta}_{A_{m,n}}(T^{m+n}, T^{-mn})\doteq
\frac{(T^{mn}-1)(T^{m+n}-1)}{(T^m-1)(T^n-1)}.$$
\end{lm}

\begin{lm}\label{lm:tor}
Let $\varphi_p : \mathbb{Z}[H_1(M)]\to \mathbb{Q}(\zeta_p)$
be a ring homomorphism such that 
$\varphi_p(T)=\zeta_p$, where $T$ is the same one in (\ref{eq:rel2}).
Then we have
$$\tau^{\varphi_p}(M)\doteq
(\zeta_p^m-1)^{-1}(\zeta_p^n-1)^{-1}.$$
\end{lm}

\noindent
{\bf Proof}\ 
By (\ref{eq:rel1}) and (\ref{eq:rel2}),
we have $\varphi_p([l_1'])\ne 1$ and $\varphi_p([l_2'])\ne 1$.
By the surgery formula II (Lemma \ref{lm:surgery2} (2)),
Theorem \ref{thm:AmnAlex}, 
Lemma \ref{lm:formula}, 
(\ref{eq:rel1}) and (\ref{eq:rel2}),
we have the result.
\qed

\section{Hyperbolicity of $b^+(m, n)$ and $b^-(m, n)$}
\label{sec:hyp}
A hyperbolic knot in $S^3$ is characterized as
a prime, non-torus and non-satellite knot.
Since the tunnel number of every Berge knot is one,
every Berge knot is prime by F.~H.~Norwood \cite{No}.
Hence we show that
both $b^+(m, n)$ and $b^-(m, n)$ are non-torus and non-satellite.
This problem is origined by 
T.~Saito and M.~Teragaito \cite{ST} for the case $b^+(m, n)$.

\medskip

To show the main theorem,
we need some known results on torus knots and satellite knots.

\begin{thm}\label{thm:Mos}
{\rm (Moser \cite{Mos})}
For an $(r, s)$-torus knot,
the surgery coefficient $p/q$
yielding a lens space satisfies $|p-qrs|=1$, and 
the resulting lens space is $L(p, -qr^2)$.
\end{thm}

\begin{thm}\label{thm:GW}
{\rm (Gordon \cite{Go};\ Wu \cite{Wu})}
A satellite knot yielding a lens space is 
a $(2, 2rs\pm 1)$-cable of an $(r, s)$-torus knot.
Then the surgery coefficient is $4rs\pm 1$, and 
the resulting lens space is $L(4rs\pm 1, \mp 4r^2)$.
\end{thm}

\begin{thm}\label{thm:Seif}
{\rm (Seifert \cite{Se})}
Let $K_{r,s}$ be an $(r, s)$-cable of a knot $K$.
Then the Alexander polynomial of $K_{r,s}$ is:
$${\mit \Delta}_{K_{r,s}}(t)\doteq
\frac{(t^{rs}-1)(t-1)}{(t^r-1)(t^s-1)}
{\mit \Delta}_K(t^r).$$
\end{thm}

\begin{lm}\label{lm:genus}
\begin{enumerate}
\item[(1)]
For an $(r, s)$-torus knot $T(r, s)$ with $r>0$ and $s>0$,
the half of the degree of the Alexander polynomial of $T(r, s)$
is $(r-1)(s-1)/2$, 
and the Reidemeister torsion of $M=(T(r, s); p)$ with $p=rs\pm 1$ is
$$\tau^{\psi_p}(M)\doteq (\zeta_p^r-1)^{-1}(\zeta_p^s-1)^{-1}$$
where a ring homomorphism 
$\psi_p : \mathbb{Z}[H_1(M)]\to \mathbb{Q}(\zeta_p)$ 
maps the meridian of $T(r, s)$ to $\zeta_p$.

\item[(2)]
For a $(2, 2rs\pm 1)$-cable of an $(r, s)$-torus knot,
denoted by $K$, the Alexander polynomial of $K$ is
$${\mit \Delta}_K(t)\doteq
\frac{(t^{2(2rs\pm 1)}-1)(t^{2rs}-1)(t-1)}
{(t^{2rs\pm 1}-1)(t^{2r}-1)(t^{2s}-1)},$$
and the Reidemeister torsion of $M=(K; p)$ with $p=4rs\pm 1$ is
$$\tau^{\psi_p}(M)\doteq (\zeta_p^{2r}-1)^{-1}(\zeta_p^{2s}-1)^{-1}$$
where a ring homomorphism 
$\psi_p : \mathbb{Z}[H_1(M)]\to \mathbb{Q}(\zeta_p)$ 
maps the meridian of $K$ to $\zeta_p$.

\end{enumerate}
\end{lm}

\noindent
{\bf Proof}\ 
(1)\ 
By Theorem \ref{thm:Seif}, we have
$${\mit \Delta}_{T(r, s)}(t)\doteq
\frac{(t^{rs}-1)(t-1)}{(t^r-1)(t^s-1)},$$
and the half of the degree of the Alexander polynomial of $T(r, s)$
is $(r-1)(s-1)/2$.
By the surgery formula II (Lemma \ref{lm:surgery2} (1)), 
we have 
$\tau^{\psi_p}(M)\doteq (\zeta_p^r-1)^{-1}(\zeta_p^s-1)^{-1}$.

\medskip

\noindent
(2)\ 
By Theorem \ref{thm:Seif}, we have the Alexander polynomial.
By the surgery formula II (Lemma \ref{lm:surgery2} (1)), 
we have 
$\tau^{\psi_p}(M)\doteq (\zeta_p^{2r}-1)^{-1}(\zeta_p^{2s}-1)^{-1}$.
\qed

\bigskip

The following is an easy but much important observation.

\begin{lm}\label{lm:odd}
For $K^{\pm}=b^+(m, n)$ and $b^-(m, n)$ respectively,
we set $p=(m+n)^2\mp mn$ and $2g=(m+n-1)^2\mp mn$.
Then $p$ is a lens surgery coefficient and $g$ is 
the half of the degree of the Alexander polynomial of $K^{\pm}$, 
and we have
$$p-2g=2(m+n)-1$$
and $p$ is odd.
\end{lm}

\noindent
{\bf Proof}\ 
By Proposition \ref{prop:genus},
Lemma \ref{lm:Amnlens} and a direct calculation,
we have the result.
\qed

\bigskip

The following is our main theorem.

\begin{thm}\label{thm:hyp}
Both $b^+(m, n)$ and $b^-(m, n)$ are hyperbolic for $2\le m<n$.
\end{thm}

\noindent
{\bf Proof}\ 
Firstly we show that both 
$K^+=b^+(m, n)$ and $K^-=b^-(m, n)$ are non-torus.
We set $p=rs\pm 1$, $g=(r-1)(s-1)/2$
and $K=K^+$ or $K^-$.
Suppose that
$K$ is an $(r, s)$-torus knot with $2\le r<s$.
Then $p=(m+n)^2-mn$ or $(m+n)^2+mn$
by Lemma \ref{lm:Amnlens} and Theorem \ref{thm:Mos}, and
the Reidemeister torsion of $M=(K; p)$ is
\begin{equation}\label{eq:tora}
\tau^{\varphi_p}(M)\doteq (\zeta_p^m-1)^{-1}(\zeta_p^n-1)^{-1}
\end{equation}
by Lemma \ref{lm:tor}, 
where a ring homomorphism 
$\varphi_p : \mathbb{Z}[H_1(M)]\to \mathbb{Q}(\zeta_p)$ 
maps the meridian of $K$ to $\zeta_p^{m+n}$ 
by (\ref{eq:rel1}).
From Lemma \ref{lm:genus} (2) and (\ref{eq:tora}),
we have
\begin{equation}\label{eq:torb}
(\zeta_p^r-1)(\zeta_p^s-1)\doteq 
(\zeta_p^{m\cdot \overline{m+n}}-1)(\zeta_p^{n\cdot \overline{m+n}}-1).
\end{equation}
From $(m+n)^2\equiv \pm mn\ (\mathrm{mod}\ \! p)$,
we have
\begin{equation}\label{eq:mn}
m\cdot \overline{m+n}\equiv \pm (m\overline{n}+1)
\quad \mbox{and}\quad
n\cdot \overline{m+n}\equiv \pm (\overline{m}n+1)\quad (\mathrm{mod}\ \! p).
\end{equation}
By applying the Franz lemma (Lemma \ref{lm:Franz}) to (\ref{eq:torb}),
we may assume
\begin{equation*}
r\equiv \pm (m\overline{n}+1)\quad \mbox{and}\quad
s\equiv \pm (\overline{m}n+1)\quad (\mathrm{mod}\ \! p).
\end{equation*}
Thus we have the following four cases:
\begin{enumerate}
\item[(i)]
$(r-1)(s-1)\equiv 1\ (\mathrm{mod}\ \! p).$

\item[(ii)]
$(r-1)(s+1)\equiv -1\ (\mathrm{mod}\ \! p).$

\item[(iii)]
$(r+1)(s-1)\equiv -1\ (\mathrm{mod}\ \! p).$

\item[(iv)]
$(r+1)(s+1)\equiv 1\ (\mathrm{mod}\ \! p).$

\end{enumerate}
We set the lefthand side of the equations above as $U$.
Then we have
$$2p-1-U\ge 2(rs-1)-1-(r+1)(s+1)=(r-1)(s-1)-5,$$
and
$$1<(r-1)(s-1)\le U<2p-1$$
except the cases $(r, s)=(2, 3)$ and $(2, 5)$ in (iv) and $p=rs-1$.

\medskip

\noindent
{\bf Case 1}\ $(r, s)\ne (2, 3)$ and $(2, 5)$.

\medskip

We have $U=p+1$ for (i) and (iv), and $U=p-1$ for (ii) and (iii).

\medskip

\noindent
(i)\ $(r-1)(s-1)=p+1$.
$$p+1-(r-1)(s-1)\ge (rs-1)+1-(r-1)(s-1)=r+s-1>0.$$
It is a contradiction.

\medskip

\noindent
(ii)\ $(r-1)(s+1)=p-1$.
$$p-1-(r-1)(s+1)\ge (rs-1)-1-(r-1)(s+1)=-r+s-1>0$$
except $s=r+1$ and $p=rs-1$.

\medskip

\noindent
(iii)\ $(r+1)(s-1)=p-1$.
$$p-1-(r+1)(s-1)\le (rs+1)-1-(r+1)(s-1)=r-s+1<0$$
except $s=r+1$ and $p=rs+1$.

\medskip

\noindent
(iv)\ $(r+1)(s+1)=p+1$.
$$p+1-(r+1)(s+1)\le (rs+1)+1-(r+1)(s+1)=-r-s+1<0.$$
It is a contradiction.

\medskip

We show the rest cases:
\begin{enumerate}
\item[(a)]
$s=r+1$ and $p=rs-1$ in (ii).

\item[(b)]
$s=r+1$ and $p=rs+1$ in (iii).

\end{enumerate}
conclude contradiction.

\medskip

\noindent
(a)\ By Proposition \ref{prop:genus}, Lemma \ref{lm:genus} (1)
and Lemma \ref{lm:odd}, we have
$$p=r(r+1)-1=(m+n)^2\mp mn\quad \mbox{and}\quad
p-2g=2r-1=2m+2n-1.$$
Thus we have
$$m+n=r\quad \mbox{and}\quad
mn=\pm(r-1).$$
Then there is no solution for $2\le m<n$.

\medskip

\noindent
(b)\ By Proposition \ref{prop:genus}, Lemma \ref{lm:genus} (1)
and Lemma \ref{lm:odd}, we have
$$p=s(s-1)+1=(m+n)^2\mp mn\quad \mbox{and}\quad
p-2g=2s-1=2m+2n-1.$$
Thus we have
$$m+n=s\quad \mbox{and}\quad
mn=\pm(s-1).$$
Then there is no solution for $2\le m<n$.

\bigskip

\noindent
{\bf Case 2}\ $(r, s)=(2, 3)$ or $(2, 5)$.

\medskip

This case corresponds to (iv).
By the assumption $2\le m<n$, we have $m+n\ge 5$.
By Proposition \ref{prop:genus}, Lemma \ref{lm:genus}
and Lemma \ref{lm:odd}, we have
$$p-2g=2m+2n-1=r+s\quad \mbox{or}\quad r+s-2.$$
If $(r, s)=(2, 3)$ or $(2, 5)$, then we have $m+n\le (r+s+1)/2\le 4$.
It is a contradiction.

\medskip

Secondly we show that both 
$K^+=b^+(m, n)$ and $K^-=b^-(m, n)$ are non-satellite.
We set $p=4rs\pm 1$, and $K=K^+$ or $K^-$.
Suppose that
$K$ is a $(2, 2rs\pm 1)$-cable of an $(r, s)$-torus knot
with $2\le r<s$.
Then $p=(m+n)^2-mn$ or $(m+n)^2+mn$
by Lemma \ref{lm:Amnlens} and Theorem \ref{thm:GW}.
In the similar way as the first case, 
by Lemma \ref{lm:genus} (2) and (\ref{eq:tora}),
we have
\begin{equation}\label{eq:torc}
(\zeta_p^{2r}-1)(\zeta_p^{2s}-1)\doteq 
(\zeta_p^{m\cdot \overline{m+n}}-1)(\zeta_p^{n\cdot \overline{m+n}}-1).
\end{equation}
By applying the Franz lemma (Lemma \ref{lm:Franz}) to (\ref{eq:torc}),
and (\ref{eq:mn}),
we may assume
\begin{equation*}
2r\equiv \pm (m\overline{n}+1)\quad \mbox{and}\quad
2s\equiv \pm (\overline{m}n+1)\quad (\mathrm{mod}\ \! p).
\end{equation*}
Thus we have the following four cases:
\begin{enumerate}
\item[(i)]
$(2r-1)(2s-1)\equiv 1\ (\mathrm{mod}\ \! p).$

\item[(ii)]
$(2r-1)(2s+1)\equiv -1\ (\mathrm{mod}\ \! p).$

\item[(iii)]
$(2r+1)(2s-1)\equiv -1\ (\mathrm{mod}\ \! p).$

\item[(iv)]
$(2r+1)(2s+1)\equiv 1\ (\mathrm{mod}\ \! p).$

\end{enumerate}
We set the lefthand side of the equations above as $U$.
Then we have
$$2p-1-U\ge 2(4rs-1)-1-(2r+1)(2s+1)
=(2r-1)(2s-1)-5\ge 10,$$
$$1<(2r-1)(2s-1)\le U<2p-1,$$
$U=p+1$ for (i) and (iv), and $U=p-1$ for (ii) and (iii).
However $U$ is odd, and $p\pm 1$ is even.
This is a contradiction.
\qed

\bigskip

Consequently from the proof of Theorem \ref{thm:hyp}, 
we have the following:

\begin{thm}\label{thm:Alex}
The Alexander polynomials of $b^+(m, n)$ and $b^-(m, n)$
for $2\le m<n$ are not those of torus knots and satellite knots
(i.e.\ The Alexander polynomials of $b^+(m, n)$ and $b^-(m, n)$
for $2\le m<n$ have already shown their hyperbolicities).
\end{thm}

We note that there may be a hyperbolic knot yielding a lens space
whose Alexander polynomial is the same as that of
a torus knot or a satellite knot.
We also note that
lens surgery coefficients of them are not always
uniquely determined by their Alexander polynomials.

\medskip

To prove Theorem \ref{thm:Alex}, we need some results.

\begin{lm}\label{lm:Kad1}
{\rm (\cite{Kad1})}
Let $K$ be a knot whose Alexander polynomial of
an $(r, s)$-torus knot with $2\le r<s$.
Suppose that $M=(K; p)$ with $p\ge 2$ has the same
Reidemeister torsions as those of a lens space.
Then we have
$$r\equiv \pm 1\quad \mbox{or}\quad
s\equiv \pm 1\quad \mbox{or}\quad
rs\equiv \pm 1\quad (\mathrm{mod}\ \! p).$$
\end{lm}

\begin{lm}\label{lm:Kad2}
{\rm (\cite{Kad2})}
Let $K$ be a knot whose Alexander polynomial of
a $(2, 2rs\pm 1)$-cable of an $(r, s)$-torus knot with $2\le r<s$.
Suppose that $M=(K; p)$ with $p\ge 2$ has the same
Reidemeister torsions as those of a lens space.
Then we have
$$r\equiv \pm 1\quad \mbox{or}\quad
s\equiv \pm 1\quad \mbox{or}\quad
2rs\equiv \pm 1\quad \mbox{or}\quad
4rs\pm 1\equiv 0\quad (\mathrm{mod}\ \! p).$$
\end{lm}

\noindent
{\bf Proof of Theorem \ref{thm:Alex}}\ 
Suppose that $K^{\pm}$ has the same Alexander polynomial 
as that of an $(r, s)$-torus knot with $2\le r<s$.
Let $p$ be a lens surgery coefficient, and
$g=(r-1)(s-1)/2$ 
the half of the degree of the Alexander polynomial of $K^{\pm}$.
By Lemma \ref{lm:Kad1}, $p=rs\pm 1$ or
\begin{equation}\label{eq:less1}
p\le \frac{rs+1}{2},
\end{equation}
or
\begin{equation}\label{eq:less2}
r=2\quad \mbox{and}\quad p=s+1.
\end{equation}
The case $p=rs\pm 1$ does not happen by the proof of Theorem \ref{thm:hyp}.
Suppose (\ref{eq:less1}).
Then we have
\begin{equation*}
2(p-2g)\le rs+1-2(r-1)(s-1)=-(r-2)(s-2)+3\le 3.
\end{equation*}
Suppose (\ref{eq:less2}).
Then we have
\begin{equation*}
p-2g=2.
\end{equation*}
By Lemma \ref{lm:odd}, we have
$$p-2g=2(m+n)-1\ge 9.$$
It is a contradiction.

\medskip

Suppose that $K^{\pm}$ has the same Alexander polynomial 
as that of a $(2, 2rs\pm 1)$-cable of an $(r, s)$-torus knot with $2\le r<s$.
Let $p$ be a lens surgery coefficient, and
$g=(r-1)(s-1)+(2rs-1\pm 1)/2$ 
the half of the degree of the Alexander polynomial of $K^{\pm}$.
By Lemma \ref{lm:Kad2}, $p=4rs\pm 1$ or
\begin{equation}\label{eq:less3}
p\le \max \left( \frac{4rs+1}{2}, 2rs+1\right)=2rs+1.
\end{equation}
The case $p=4rs\pm 1$ does not happen by the proof of Theorem \ref{thm:hyp}.
We suppose (\ref{eq:less3}).
Then we have
\begin{equation*}
p-2g\le 2rs+1-2(r-1)(s-1)-(2rs-2)=-2(r-1)(s-1)+3\le -1.
\end{equation*}
By Lemma \ref{lm:odd}, we have a contradiction.
\qed

\section{Identification of $b^+(m, n)$ and $b^-(m, n)$}
\label{sec:detect}
We ask what kind of information identify a Berge knot of type VII or VIII.
We consider the cases that 
(i) the resulting lens space up to orientations 
by an odd integer surgery along the knot, 
(ii) a pair of an odd integer surgery $p$ and 
the half of the degree of the Alexander polynomial $g$ of the knot, and
(iii) one of $p$ and $g$.

\medskip

A background of the case (i) is the following fact:

\begin{thm}\label{thm:ST}
{\rm (Saito and Teragaito \cite[Theorem 1.1]{ST})}
There are infinite pairs of two distinct knots in $S^3$ 
whose same integer surgeries yield homeomorphic lens spaces.
Then we can take the two knots from the classes of
torus knots, satellite knots, and hyperbolic knots arbitrarily
except the case that both knots are satellite knots.
\end{thm}

The Cyclic Surgery Theorem \cite{CGLS} says that
the number of lens surgery coefficients of a hyperbolic knot is 
at most $2$, and they are successing integers if it is $2$.
Hence the odd lens surgery coefficient
of a hyperbolic knot is uniquely determined if it exists.

\begin{lm}\label{lm:lens}
The resulting lens space up to orientations 
by an odd integer surgery along
a hyperbolic Berge knot of type VII or VIII identifies uniquely 
the standard parameter.
\end{lm}

\noindent
{\bf Proof}\ 
For $b^+(m, n)$ or $b^-(m, n)$ with $2\le m<n$
(i.e.\ a hyperbolic Berge knot of type VII or VIII),
$p=(m+n)^2\mp mn$ is the odd lens surgery coefficient.
Then the result of $p$-surgery is $L(p, m\overline{n})$
by Lemma \ref{lm:Amnlens}.
We set its standard parameter as $(\varepsilon, m, n)$.
Suppose that there is a hyperbolic Berge knot of type VII or VIII
yielding a homeomorphic lens space
with its standard parameter $(\varepsilon', m', n')$.
Then we have
$$m\overline{n}\equiv \pm m'\overline{n'}\quad
\mbox{or}\quad \pm \overline{m'}n'\quad
(\mathrm{mod}\ \! p),$$
and
$$mn'+m'n\quad \mbox{or}\quad 
mn'-m'n\quad \mbox{or}\quad
mm'+nn'\quad \mbox{or}\quad 
mm'-nn'\equiv 0\quad (\mathrm{mod}\ \! p).$$
By the Cauchy-Schwartz inequation, we have
\begin{eqnarray*}
0 & \le &
\min (mn'+m'n, |mn'-m'n|, mm'+nn', |mm'-nn'|)\\
& < &
\max (mn'+m'n, |mn'-m'n|, mm'+nn', |mm'-nn'|)\\
& \le & \sqrt{m^2+n^2}\cdot \sqrt{(m')^2+(n')^2}\\
& \le & \max (m^2+n^2, (m')^2+(n')^2)<p,
\end{eqnarray*}
and hence only the case $|mn'-m'n|=0$ is possible.
This is equivalent to $(m, n)=(m', n')$.
It deduces $\varepsilon=\varepsilon'$.
\qed

\bigskip

By Lemma \ref{lm:lens}, we have the following.

\begin{co}\label{co:mirror}
The standard parameter with $2\le m<n$ identifies uniquely 
a hyperbolic Berge knot of type VII or VIII up to mirror images.
\end{co}

For the case (ii), we set
$$
\begin{array}{ccl}
\mathcal{S} & = & 
\{(\varepsilon, m, n)\in \{ +, -\}\times \mathbb{N}\times \mathbb{N}\ 
|\ 1\le m\le n,\ \gcd(m, n)=1\},\\
\mathcal{S}_2 & = & 
\{(\varepsilon, m, n)\in \{ +, -\}\times \mathbb{N}\times \mathbb{N}\ 
|\ 2\le m<n,\ \gcd(m, n)=1\},
\medskip \\
\mathcal{B}_{+,1} & = & \{b^+(1, n)\ |\ 1\le n\}
=\{T(n, n+1)\ |\ 1\le n\},\\
\mathcal{B}_{-,1} & = & \{b^-(1, n)\ |\ 1\le n\}
=\{T(n+1, n+2)\ |\ 1\le n\},\medskip \\
\mathcal{B}_2 & = & \{b^+(m, n)\ |\ 2\le m<n,\ \gcd(m, n)=1\}\\
 &  & \cup \{b^-(m, n)\ |\ 2\le m<n,\ \gcd(m, n)=1\},\medskip \\
\mathcal{B} & = & \mathcal{B}_{+,1}\cup \mathcal{B}_{-,1}\cup \mathcal{B}_2,
\end{array}
$$
where $\mathcal{S}$ is the set of the standard parameters, 
and knots are considered up to orientations and mirror images.
We define maps
$$
F : \mathcal{S}\to \mathbb{N}\times \mathbb{N}
\quad \mbox{and} \quad
G : \mathcal{S}\to \mathcal{B}
$$
by 
\begin{eqnarray*}
F(\varepsilon, m, n) & = & 
(p, 2g)=((m+n)^2-\varepsilon mn, (m+n-1)^2-\varepsilon mn),\\
G(\varepsilon, m, n) & = & 
b^+(m, n)\ \mbox{for $\varepsilon=+$}\ \mbox{and}\ 
b^-(m, n)\ \mbox{for $\varepsilon=-$}.
\end{eqnarray*}

\begin{thm}\label{thm:detect}
\begin{enumerate}
\item[(1)]
The map $F$ is injective.

\item[(2)]
$\mathcal{B}_{+,1}\supset \mathcal{B}_{-,1}$ and
$\mathcal{B}_{+,1}\setminus \mathcal{B}_{-,1}
=\{b^+(1, 1)\}$.

\item[(3)]
The map $G$ is surjective, and
for any $K\in \mathcal{B}$,
$\sharp \{G^{-1}(K)\}=1$ or $2$
where $\sharp \{ \cdot \}$ implies the cardinality.
$\sharp \{G^{-1}(K)\}=2$ if and only if
$K\in \mathcal{B}_{-,1}$
(i.e.\ For some $n\ge 2$, 
$K=T(n, n+1)$ and
$G^{-1}(K)=\{ (+, 1, n), (-, 1, n-1)\}$).

\end{enumerate}
\end{thm}

\noindent
{\bf Proof}\ 
(1)\ By Lemma \ref{lm:odd}, 
we have
$$p-2g=2(m+n)-1
\quad \mbox{and} \quad
p-(m+n)^2=-\varepsilon mn.$$
Since $m+n$ and $mn$ are uniquely determined from 
$\varepsilon$ and $(p, 2g)$, 
$m$ and $n$ are also uniquely determined from 
$\varepsilon$ and $(p, 2g)$.

\medskip

Suppose that $F(+, m, n)=F(-, m', n')=(p, 2g)$.
By the relations above, we have $m+n=m'+n'$ and
$$p-(m+n)^2=-mn=m'n'.$$
It is a contradiction.
Therefore $F$ is injective.

\bigskip

\noindent
(2)\ Since $b^+(1, n)=b^-(1, n-1)=T(n, n+1)$,
we have the result.

\bigskip

\noindent
(3)\ Surjectivity of $G$ is clear.
By Theorem \ref{thm:hyp} and the Cyclic Surgery Theorem \cite{CGLS},
the odd lens surgery coefficient of an element in 
$\mathcal{B}_2$ is uniquely determined as $p$.
By (1) and (2), we have the result.
\qed

\bigskip

By Theorem \ref{thm:detect}, we have the following.

\begin{co}\label{co:pg}
The map $G\circ F^{-1} : F(\mathcal{S}_2)\to \mathcal{B}$
is injective, and its image is $\mathcal{B}_2$.
That is, a pair $(p, 2g)$ with $2\le m<n$ identifies uniquely 
a hyperbolic Berge knot of type VII or VIII up to mirror images.
\end{co}

\begin{rem}\label{rem:Yam}
{\rm
(1)\ Lemma \ref{lm:lens} and Corollary \ref{co:mirror}
say that $p$ and $\pm m\overline{n}$\ $(\mathrm{mod}\ \! p)$
with $2\le m<n$ determine $(\varepsilon, m, n)$ and 
a hyperbolic Berge knot of type VII or VIII uniquely.
Theorem \ref{thm:detect} and Corollary \ref{co:pg} say that
$p$ and $g$ with $2\le m<n$ determine $(\varepsilon, m, n)$  and
a hyperbolic Berge knot of type VII or VIII uniquely.

\medskip

\noindent
(2)
Y.~Yamada \cite{Yam2} defines $k^+(m, n)(=b^+(m, n))$ and 
$k^-(m, m+n)(=b^-(m, n))$ for any coprime pair $(m, n)$, and
their mirror images $l^+(m, n)$ and $l^-(m, m+n)$, respectively.
They have unique reprsentatives in $\mathcal{B}$
up to orientations and mirror images
(see Lemma \ref{lm:primitive} (3)).
Since a double torus knot is strongly invertible, 
a non-trivial torus knot is non-amphicheiral, and
a non-torus knot yielding a lens space is non-amphicheiral
by the Cyclic Surgery Theorem \cite{CGLS}, 
the extended class can be also classified completely.}
\end{rem}

For the case (iii), we define maps
$$
\mathrm{proj}_1 : \mathbb{N}\times \mathbb{N}\to \mathbb{N}
\quad \mbox{and}\quad
\mathrm{proj}_2 : \mathbb{N}\times \mathbb{N}\to \mathbb{N}
$$
by 
$$\mathrm{proj}_1(x, y)=x
\quad \mbox{and}\quad
\mathrm{proj}_2(x, y)=y,$$
respectively.
We set
$H_{\varepsilon,i}=\mathrm{proj}_i\circ F_{\varepsilon}$
and
$H_i=\mathrm{proj}_i\circ F$\ 
$(\varepsilon=+, -;\ i=1, 2)$.
Then we study
$\mathrm{Im}\ \! (H_{\varepsilon,i})$,
$\mathrm{Im}\ \! (H_{+,i})\cap \mathrm{Im}\ \! (H_{-,i})$
and $(H_2)^{-1}(2g)$.

\medskip

We need number theoretical results on
$\mathbb{Q}(\sqrt{-3})$ and $\mathbb{Q}(\sqrt{5})$.
Let $O_+$ and $O_-$ be the integer rings of 
$\mathbb{Q}(\sqrt{-3})$ and $\mathbb{Q}(\sqrt{5})$
respectively, and
$U_+=O_+^{\times}$ and $U_-=O_-^{\times}$ 
the unit groups of $O_+$ and $O_-$ respectively.
We set
$$\omega_+=\frac{1+\sqrt{-3}}{2},\quad
\overline{\omega_+}=\frac{1-\sqrt{-3}}{2},\quad
\omega_-=\frac{1+\sqrt{5}}{2}\quad \mbox{and}\quad
\overline{\omega_-}=\frac{1-\sqrt{5}}{2}.$$
Then $\mathbb{Q}(\sqrt{-3})=\mathbb{Q}(\omega_+)$ 
and $\mathbb{Q}(\sqrt{5})=\mathbb{Q}(\omega_-)$, and
their Galois groups are
$$\mathrm{Gal}\ \! (\mathbb{Q}(\omega_+)/\mathbb{Q})
=\langle \sigma_+\ |\ \sigma_+^2=1\rangle
\cong \mathbb{Z}/2\mathbb{Z}
\quad \mbox{and}\quad 
\mathrm{Gal}\ \! (\mathbb{Q}(\omega_-)/\mathbb{Q})
=\langle \sigma_-\ |\ \sigma_-^2=1\rangle
\cong \mathbb{Z}/2\mathbb{Z}$$
where $\sigma_{\pm}(\omega_{\pm})=\overline{\omega_{\pm}}$.
For $x\in \mathbb{Q}(\omega_{\pm})$,
we set $\overline{x}=\sigma_{\pm}(x)$.
It is easy to see that
$$\omega_++\overline{\omega_+}=1,\ 
\omega_+\overline{\omega_+}=1,\ 
(\omega_+)^2=-\overline{\omega_+},\ 
(\omega_+)^3=-1,\ 
(\omega_+)^6=1,$$
and
$$\omega_-+\overline{\omega_-}=1,\ 
\omega_-\overline{\omega_-}=-1,\ 
(\omega_-)^2=1+\omega_-=\frac{3+\sqrt{5}}{2}.$$
Let $\{ a_k\}_{k\in \mathbb{Z}}$ be the {\it Fibonacci series}
defined by
$$a_{k+2}=a_{k+1}+a_k\quad \mbox{and}\quad 
a_1=a_2=1.$$

\medskip

The following is elementary results on a theory of quadratic fields
except (3).

\begin{prop}\label{prop:quad}
\begin{enumerate}
\item[(1)]
$O_+=\mathbb{Z}[\omega_+]$ and
$O_-=\mathbb{Z}[\omega_-]$ as rings, and
$$O_+\cong \mathbb{Z}\oplus \mathbb{Z}\cdot \omega_+
\quad \mbox{and}\quad 
O_-\cong \mathbb{Z}\oplus \mathbb{Z}\cdot \omega_-
\cong \mathbb{Z}\oplus \mathbb{Z}\cdot (\omega_-)^2$$
as abelian groups.
Both $O_+$ and $O_-$ have the class number one
(i.e.\ principal ideal domains).

\item[(2)]
$U_+=\langle \omega_+\ |\ (\omega_+)^6=1\rangle
\cong \mathbb{Z}/6\mathbb{Z}$
and
$U_-=\langle \omega_-, -1\ |\ (-1)^2=1\rangle
\cong \mathbb{Z}\oplus \mathbb{Z}/2\mathbb{Z}.$

\item[(3)]
\begin{enumerate}
\item[(a)]
For $k\ge 1$, $a_k$ is a positive integer.

\item[(b)]
For $k\in \mathbb{Z}$, 
$$a_k=\frac{(\omega_-)^k-(\overline{\omega_-})^k}
{\omega_--\overline{\omega_-}},\quad
a_{-k}=(-1)^{k+1}a_k,$$
and
$$(\omega_-)^k=a_{k-1}+a_k\omega_-.$$

\end{enumerate}

\item[(4)]
Let $\ell$ be a positive prime of $\mathbb{Z}$.
\begin{enumerate}
\item[(a)]
$\ell$ is not a prime in $O_+$ if and only if
$\ell=3$ or $\equiv 1\ (\mathrm{mod}\ \! 3)$.
Their prime factorizations in $O_+$ are
$$3=-(\sqrt{-3})^2\quad \mbox{or}\quad 
\ell=\frak{l}\cdot \overline{\frak{l}}\quad 
(\frak{l}\in O_+, \frak{l}\ne \overline{\frak{l}}),$$
respectively.

\item[(b)]
$\ell$ is not a prime in $O_-$ if and only if
$\ell=5$ or $\equiv \pm 1\ (\mathrm{mod}\ \! 5)$.
Their prime factorizations in $O_-$ are
$$5=(\sqrt{5})^2\quad \mbox{or}\quad 
\ell=\frak{l}\cdot \overline{\frak{l}}\quad 
(\frak{l}\in O_-, \frak{l}\ne \overline{\frak{l}}),$$
respectively.

\end{enumerate}
\end{enumerate}
\end{prop}

We say that an element 
$b+c\omega_{\pm}\in O_{\pm}$\ $(b, c\in \mathbb{Z})$
is {\it primitive} if $\gcd(b, c)=1$.

\begin{lm}\label{lm:primitive}
\begin{enumerate}
\item[(1)]
For $x=b+c\omega_{\pm}\in O_{\pm}$\ $(b, c\in \mathbb{Z})$, we set
$x_k=x(\omega_{\pm})^{k-1}$\ $(k\in \mathbb{Z})$.
\begin{enumerate}
\item[(a)]
For $x=b+c\omega_+\in O_+$\ $(b, c\in \mathbb{Z})$, we have
$$\begin{array}{rclrcl}
x=x_1 & = & b+c\omega_+,\quad & 
x_2 & = & -c+(b+c)\omega_+, \\
x_3 & = & -(b+c)+b\omega_+,\quad &
x_4 & = & -b-c\omega_+, \\
x_5 & = & c-(b+c)\omega_+,\quad &
x_6 & = & (b+c)-b\omega_+, \\
\overline{x}=\overline{x_1} & = & (b+c)-c\omega_+,\quad & 
\overline{x_2} & = & b-(b+c)\omega_+, \\
\overline{x_3} & = & -c-b\omega_+,\quad &
\overline{x_4} & = & -(b+c)+c\omega_+, \\
\overline{x_5} & = & -b+(b+c)\omega_+,\quad &
\overline{x_6} & = & c+b\omega_+,
\end{array}$$
and $x_{k+6}=x_k$.

\item[(b)]
For $x=b+c\omega_-\in O_-$\ $(b, c\in \mathbb{Z})$, we set
$x_k=b_k+b_{k+1}\omega_-$\ $(b_k\in \mathbb{Z})$.
Then $\{ b_k\}_{k\in \mathbb{Z}}$ is the Fibonacci series
such that
$$b_{k+2}=b_{k+1}+b_k,\quad b_1=b\quad \mbox{and}\quad b_2=c,$$
and we have
$b_k=ba_{k-2}+ca_{k-1}$
where $\{ a_k\}_{k\in \mathbb{Z}}$ is the same as in
Proposition \ref{prop:quad} (3), and
$\overline{x_k}=b_{k+2}-b_{k+1}\omega_-$.

\end{enumerate}

\item[(2)]
For $x=b+c\omega_{\pm}\in O_{\pm}$\ $(b, c\in \mathbb{Z})$, we set
$$\mathrm{Orb}\ \! (x)=\{xu,\ \overline{x}u\ |\ u\in U_{\pm}\}.$$
Then one element of $\mathrm{Orb}\ \! (x)$ is primitive
if and only if every element of $\mathrm{Orb}\ \! (x)$ is primitive.

\item[(3)]
Let $x\in O_{\pm}$ be a primitive element.
Then we have the following.
\begin{enumerate}
\item[(a)]
There is a unique element 
$b+c\omega_+\in \mathrm{Orb}\ \! (x)$\ $(b, c\in \mathbb{Z})$
such that $1\le b\le c$ (and $\gcd(b, c)=1$).

\item[(b)]
There is a unique element 
$b+c(\omega_-)^2
=(b+c)+c\omega_-\in \mathrm{Orb}\ \! (x)$\ $(b, c\in \mathbb{Z})$
such that $1\le b\le c$ (and $\gcd(b, c)=1$).

\end{enumerate}

\item[(4)]
Let $x, y\in O_{\pm}$ be two primitive elements.
Then $xy$ is also primitive if and only if
for every prime factor $\frak{l}$ of $x$, 
$\overline{\frak{l}}$ is not a prime factor of $y$.

\end{enumerate}
\end{lm}

\noindent
{\bf Proof}\ 
(1)\ We can show by straight calculations.

\bigskip

\noindent
(2)\ By (1) and Proposition \ref{prop:quad} (2), we have the result.

\bigskip

\noindent
(3)\ By (1) and (2), they can be checked without difficulty.

\bigskip

\noindent
(4)\ By Proposition \ref{prop:quad} (4), we have the result.
\qed

\bigskip

We can factorize
\begin{eqnarray}\label{eq:factor}
\begin{matrix}
m^2+mn+n^2 & = & 
(m+n\omega_+)(m+n\overline{\omega_+}), \hfill \medskip \\
m^2+3mn+n^2 & = & 
\{ m+n(\omega_-)^2\} \{m+n(\overline{\omega_-})^2\}\hfill \medskip \\
& = & \{(m+n)+n\omega_-\}\{(m+n)+n\overline{\omega_-}\}\hfill
\end{matrix}
\end{eqnarray}
in $O_+$ and $O_-$, respectively.
Note that every factor in the righthand side of (\ref{eq:factor})
is primitive, and satisfies the conditions of 
Lemma \ref{lm:primitive} (3).

\begin{lm}\label{lm:prime}
\begin{enumerate}
\item[(1)]
{\rm (Berge \cite[Theorem 4]{Ber})}
Let $p$ be a positive integer with a prime factorization
$p=3^e\ell_1^{e_1}\cdots \ell_r^{e_r}$
where $\ell_1, \ldots, \ell_r$ are distinct primes
other than $3$, $e\ge 0$ and $e_i\ge 1$\ $(i=1, \ldots, r)$.
Then $p\in \mathrm{Im}\ \! (H_{+,1})$ if and only if
$e=0$ or $1$, and every $\ell_i\equiv 1\ (\mathrm{mod}\ \! 3)$.
Moreover the number of elements of the set
$(H_{+,1})^{-1}(p)$ is $2^{r-1}$.

\item[(2)]
{\rm (Berge \cite[Theorem 5]{Ber})}
Let $p$ be a positive integer with a prime factorization
$p=5^e\ell_1^{e_1}\cdots \ell_r^{e_r}$
where $\ell_1, \ldots, \ell_r$ are distinct primes
other than $5$, $e\ge 0$ and $e_i\ge 1$\ $(i=1, \ldots, r)$.
Then $p\in \mathrm{Im}\ \! (H_{-,1})$ if and only if
$e=0$ or $1$, and every 
$\ell_i\equiv 1$ or $-1\ (\mathrm{mod}\ \! 5)$.
Moreover the number of elements of the set
$(H_{-,1})^{-1}(p)$ is $2^{r-1}$.

\item[(3)]
Let $p$ be a positive integer with a prime factorization
$p=\ell_1^{e_1}\cdots \ell_r^{e_r}$
where $\ell_1, \ldots, \ell_r$ are distinct primes and $e_i\ge 1$\ $(i=1, \ldots, r)$.
Then $p\in \mathrm{Im}\ \! (H_{+,1})\cap \mathrm{Im}\ \! (H_{-,1})$ 
if and only if every $\ell_i\equiv 1$ or $4\ (\mathrm{mod}\ \! 15)$.

\end{enumerate}
\end{lm}

\noindent
{\bf Proof}\ 
(1)\ By Proposition \ref{prop:quad} (4) (a), 
Lemma \ref{lm:primitive} (4),
and (\ref{eq:factor}), we have the result.

\bigskip

\noindent
(2)\ By Proposition \ref{prop:quad} (4) (b), 
Lemma \ref{lm:primitive} (4),
and (\ref{eq:factor}), we have the result.

\bigskip

\noindent
(3)\ By (1) and (2), we have the result.
\qed

\bigskip

By Proposition \ref{prop:genus}, Lemma \ref{lm:Amnlens} and
Lemma \ref{lm:prime},
we make a complete table of $p$, $(\varepsilon, m, n)$ and $g$ for $p\le 500$.

\newpage

\begin{minipage}[t]{7cm}
\begin{center}
\begin{tabular}{|c@{\quad\vrule width0.8pt\quad}c|c|}
\hline
$p$ & $(\varepsilon, m, n)$ & $g$\\
\noalign{\hrule height 0.8pt}
$7$ & $(+, 1, 2)$ & $1$\\
\hline
$11$ & $(-, 1, 2)$ & $3$\\
\hline
$13$ & $(+, 1, 3)$ & $3$\\
\hline
$19$ & $(+, 2, 3)$ & $5$\\
\hline
 & $(-, 1, 3)$ & $6$\\
\hline
$21$ & $(+, 1, 4)$ & $6$\\
\hline
$29$ & $(-, 1, 4)$ & $10$\\
\hline
$31$ & $(+, 1, 5)$ & $10$\\
\hline
 & $(-, 2, 3)$ & $11$\\
\hline
$37$ & $(+, 3, 4)$ & $12$\\
\hline
$39$ & $(+, 2, 5)$ & $13$\\
\hline
$41$ & $(-, 1, 5)$ & $15$\\
\hline
$43$ & $(+, 1, 6)$ & $15$\\
\hline
$49$ & $(+, 3, 5)$ & $17$\\
\hline
$55$ & $(-, 1, 6)$ & $21$\\
\hline
$57$ & $(+, 1, 7)$ & $21$\\
\hline
$59$ & $(-, 2, 5)$ & $23$\\
\hline
$61$ & $(+, 4, 5)$ & $22$\\
\hline
 & $(-, 3, 4)$ & $24$\\
\hline
$67$ & $(+, 2, 7)$ & $25$\\
\hline
$71$ & $(-, 1, 7)$ & $28$\\
\hline
$73$ & $(+, 1, 8)$ & $27$\\
\hline
$79$ & $(+, 3, 7)$ & $30$\\
\hline
 & $(-, 3, 5)$ & $32$\\
\hline
$89$ & $(-, 1, 8)$ & $36$\\
\hline
$91$ & $(+, 1, 9)$ & $36$\\
\hline
 & $(+, 5, 6)$ & $35$\\
\hline
$93$ & $(+, 4, 7)$ & $36$\\
\hline
$95$ & $(-, 2, 7)$ & $39$\\
\hline
$97$ & $(+, 3, 8)$ & $38$\\
\hline
$101$ & $(-, 4, 5)$ & $42$\\
\hline
$103$ & $(+, 2, 9)$ & $41$\\
\hline
$109$ & $(+, 5, 7)$ & $43$\\
\hline
 & $(-, 1, 9)$ & $45$\\
\hline
$111$ & $(+, 1, 10)$ & $45$\\
\hline
$121$ & $(-, 3, 7)$ & $51$\\
\hline
$127$ & $(+, 6, 7)$ & $51$\\
\hline
$129$ & $(+, 5, 8)$ & $52$\\
\hline
$131$ & $(-, 1, 10)$ & $55$\\
\hline
$133$ & $(+, 1, 11)$ & $55$\\
\hline
\end{tabular}
\end{center}
\end{minipage}
\ 
\begin{minipage}[t]{7cm}
\begin{center}
\begin{tabular}{|c@{\quad\vrule width0.8pt\quad}c|c|}
\hline
$p$ & $(\varepsilon, m, n)$ & $g$\\
\noalign{\hrule height 0.8pt}
 & $(+, 4, 9)$ & $54$\\
\hline
$139$ & $(+, 3, 10)$ & $57$\\
\hline
 & $(-, 2, 9)$ & $59$\\
\hline
$145$ & $(-, 3, 8)$ & $62$\\
\hline
$147$ & $(+, 2, 11)$ & $61$\\
\hline
$149$ & $(-, 4, 7)$ & $64$\\
\hline
$151$ & $(+, 5, 9)$ & $62$\\
\hline
 & $(-, 5, 6)$ & $65$\\
\hline
$155$ & $(-, 1, 11)$ & $66$\\
\hline
$157$ & $(+, 1, 12)$ & $66$\\
\hline
$163$ & $(+, 3, 11)$ & $68$\\
\hline
$169$ & $(+, 7, 8)$ & $70$\\
\hline
$179$ & $(-, 5, 7)$ & $78$\\
\hline
$181$ & $(+, 4, 11)$ & $76$\\
\hline
 & $(-, 1, 12)$ & $78$\\
\hline
$183$ & $(+, 1, 13)$ & $78$\\
\hline
$191$ & $(-, 2, 11)$ & $83$\\
\hline
$193$ & $(+, 7, 9)$ & $81$\\
\hline
$199$ & $(+, 2, 13)$ & $85$\\
\hline
 & $(-, 3, 10)$ & $87$\\
\hline
$201$ & $(+, 5, 11)$ & $85$\\
\hline
$205$ & $(-, 4, 9)$ & $90$\\
\hline
$209$ & $(-, 1, 13)$ & $91$\\
\hline
 & $(-, 5, 8)$ & $92$\\
\hline
$211$ & $(+, 1, 14)$ & $91$\\
\hline
 & $(-, 6, 7)$ & $93$\\
\hline
$217$ & $(+, 3, 13)$ & $93$\\
\hline
 & $(+, 8, 9)$ & $92$\\
\hline
$219$ & $(+, 7, 10)$ & $93$\\
\hline
$223$ & $(+, 6, 11)$ & $95$\\
\hline
$229$ & $(+, 5, 12)$ & $98$\\
\hline
 & $(-, 3, 11)$ & $101$\\
\hline
$237$ & $(+, 4, 13)$ & $102$\\
\hline
$239$ & $(-, 1, 14)$ & $105$\\
\hline
$241$ & $(+, 1, 15)$ & $105$\\
\hline
 & $(-, 5, 9)$ & $107$\\
\hline
$247$ & $(+, 3, 14)$ & $107$\\
\hline
 & $(+, 7, 11)$ & $106$\\
\hline
$251$ & $(-, 2, 13)$ & $111$\\
\hline
$259$ & $(+, 2, 15)$ & $113$\\
\hline
\end{tabular}
\end{center}
\end{minipage}

\newpage

\begin{minipage}[t]{7cm}
\begin{center}
\begin{tabular}{|c@{\quad\vrule width0.8pt\quad}c|c|}
\hline
$p$ & $(\varepsilon, m, n)$ & $g$\\
\noalign{\hrule height 0.8pt}
 & $(+, 5, 13)$ & $112$\\
\hline
$269$ & $(-, 4, 11)$ & $120$\\
\hline
$271$ & $(+, 9, 10)$ & $117$\\
\hline
 & $(-, 1, 15)$ & $120$\\
\hline
$273$ & $(+, 1, 16)$ & $120$\\
\hline
 & $(+, 8, 11)$ & $118$\\
\hline
$277$ & $(+, 7, 12)$ & $120$\\
\hline
$281$ & $(-, 7, 8)$ & $126$\\
\hline
$283$ & $(+, 6, 13)$ & $123$\\
\hline
$291$ & $(+, 5, 14)$ & $127$\\
\hline
$295$ & $(-, 3, 13)$ & $132$\\
\hline
$301$ & $(+, 4, 15)$ & $132$\\
\hline
 & $(+, 9, 11)$ & $131$\\
\hline
$305$ & $(-, 1, 16)$ & $136$\\
\hline
$307$ & $(+, 1, 17)$ & $136$\\
\hline
$309$ & $(+, 7, 13)$ & $135$\\
\hline
$311$ & $(-, 5, 11)$ & $140$\\
\hline
$313$ & $(+, 3, 16)$ & $138$\\
\hline
$319$ & $(-, 2, 15)$ & $143$\\
\hline
 & $(-, 7, 9)$ & $144$\\
\hline
$327$ & $(+, 2, 17)$ & $145$\\
\hline
$331$ & $(+, 10, 11)$ & $145$\\
\hline
 & $(-, 3, 14)$ & $149$\\
\hline
$337$ & $(+, 8, 13)$ & $148$\\
\hline
$341$ & $(-, 1, 17)$ & $153$\\
\hline
 & $(-, 4, 13)$ & $154$\\
\hline
$343$ & $(+, 1, 18)$ & $153$\\
\hline
$349$ & $(+, 3, 17)$ & $155$\\
\hline
 & $(-, 5, 12)$ & $158$\\
\hline
$355$ & $(-, 6, 11)$ & $161$\\
\hline
$359$ & $(-, 2, 15)$ & $163$\\
\hline
$361$ & $(+, 5, 16)$ & $160$\\
\hline
 & $(-, 8, 9)$ & $164$\\
\hline
$367$ & $(+, 9, 13)$ & $162$\\
\hline
$373$ & $(+, 4, 17)$ & $166$\\
\hline
$379$ & $(+, 7, 15)$ & $168$\\
\hline
 & $(-, 1, 18)$ & $171$\\
\hline
$381$ & $(+, 1, 19)$ & $171$\\
\hline
$389$ & $(-, 5, 13)$ & $177$\\
\hline
$395$ & $(-, 2, 17)$ & $179$\\
\hline
\end{tabular}
\end{center}
\end{minipage}
\ 
\begin{minipage}[t]{7cm}
\begin{center}
\begin{tabular}{|c@{\quad\vrule width0.8pt\quad}c|c|}
\hline
$p$ & $(\varepsilon, m, n)$ & $g$\\
\noalign{\hrule height 0.8pt}
$397$ & $(+, 11, 12)$ & $176$\\
\hline
$399$ & $(+, 5, 17)$ & $178$\\
\hline
 & $(+, 10, 13)$ & $177$\\
\hline
$401$ & $(-, 7, 11)$ & $183$\\
\hline
$403$ & $(+, 2, 19)$ & $181$\\
\hline
 & $(+, 9, 14)$ & $179$\\
\hline
$409$ & $(+, 8, 15)$ & $182$\\
\hline
 & $(-, 3, 16)$ & $186$\\
\hline
$417$ & $(+, 7, 16)$ & $186$\\
\hline
$419$ & $(-, 1, 19)$ & $190$\\
\hline
$421$ & $(+, 1, 20)$ & $190$\\
\hline
 & $(-, 4, 15)$ & $192$\\
\hline
$427$ & $(+, 3, 19)$ & $192$\\
\hline
 & $(+, 6, 17)$ & $191$\\
\hline
$431$ & $(-, 5, 14)$ & $197$\\
\hline
$433$ & $(+, 11, 13)$ & $193$\\
\hline
$439$ & $(+, 5, 18)$ & $197$\\
\hline
 & $(-, 6, 13)$ & $201$\\
\hline
$449$ & $(-, 8, 11)$ & $206$\\
\hline
$451$ & $(-, 3, 17)$ & $206$\\
\hline
 & $(-, 9, 10)$ & $207$\\
\hline
$453$ & $(+, 4, 19)$ & $204$\\
\hline
$455$ & $(-, 7, 12)$ & $204$\\
\hline
$457$ & $(+, 7, 17)$ & $205$\\
\hline
$461$ & $(-, 1, 20)$ & $210$\\
\hline
$463$ & $(+, 1, 21)$ & $210$\\
\hline
$469$ & $(+, 3, 20)$ & $212$\\
\hline
 & $(+, 12, 13)$ & $210$\\
\hline
$471$ & $(+, 11, 14)$ & $211$\\
\hline
$479$ & $(-, 2, 19)$ & $219$\\
\hline
$481$ & $(+, 5, 19)$ & $217$\\
\hline
 & $(+, 9, 16)$ & $216$\\
\hline
$487$ & $(+, 2, 21)$ & $221$\\
\hline
$489$ & $(+, 8, 17)$ & $220$\\
\hline
$491$ & $(-, 7, 13)$ & $226$\\
\hline
$499$ & $(+, 7, 18)$ & $225$\\
\hline
 & $(-, 9, 11)$ & $230$\\
\hline
\end{tabular}
\end{center}
\end{minipage}

\newpage

From the table, we raise examples of $g$ such that
$\sharp \{(H_2)^{-1}(2g)\} \ge 2$,
where $\sharp \{ \cdot \}$ implies the cardinality,
except the trivial case 
$(H_2)^{-1}(n(n-1))=\{ (+, 1, n), (-, 1, n-1)\}$
(Then $n\ge 2$, $H_1(+, 1, n)=n^2+n+1$ and $H_1(-, 1, n-1)=n^2+n-1$).

\bigskip

\begin{minipage}[t]{7cm}
\begin{center}
\begin{tabular}{|c@{\quad\vrule width0.8pt\quad}c|c|}
\hline
$p$ & $(\varepsilon, m, n)$ & $g$\\
\noalign{\hrule height 0.8pt}
$89$ & $(-, 1, 8)$ & $36$\\
\hline
$91$ & $(+, 1, 9)$ & $36$\\
\hline
$93$ & $(+, 4, 7)$ & $36$\\
\hline
\hline
$121$ & $(-, 3, 7)$ & $51$\\
\hline
$127$ & $(+, 6, 7)$ & $51$\\
\hline
\hline
$145$ & $(-, 3, 8)$ & $62$\\
\hline
$151$ & $(+, 5, 9)$ & $62$\\
\hline
\hline
$179$ & $(-, 5, 7)$ & $78$\\
\hline
$181$ & $(-, 1, 12)$ & $78$\\
\hline
$183$ & $(+, 1, 13)$ & $78$\\
\hline
\hline
$199$ & $(+, 2, 13)$ & $85$\\
\hline
$201$ & $(+, 5, 11)$ & $85$\\
\hline
\hline
$209$ & $(-, 5, 8)$ & $92$\\
\hline
$217$ & $(+, 8, 9)$ & $92$\\
\hline
\hline
$211$ & $(-, 6, 7)$ & $93$\\
\hline
$217$ & $(+, 3, 13)$ & $93$\\
\hline
$219$ & $(+, 7, 10)$ & $93$\\
\hline
\hline
$241$ & $(-, 5, 9)$ & $107$\\
\hline
$247$ & $(+, 3, 14)$ & $107$\\
\hline
\hline
$269$ & $(-, 4, 11)$ & $120$\\
\hline
$271$ & $(-, 1, 15)$ & $120$\\
\hline
$273$ & $(+, 1, 16)$ & $120$\\
\hline
$277$ & $(+, 7, 12)$ & $120$\\
\hline
\end{tabular}
\end{center}
\end{minipage}
\quad
\begin{minipage}[t]{7cm}
\begin{center}
\begin{tabular}{|c@{\quad\vrule width0.8pt\quad}c|c|}
\hline
$p$ & $(\varepsilon, m, n)$ & $g$\\
\noalign{\hrule height 0.8pt}
$295$ & $(-, 3, 13)$ & $132$\\
\hline
$301$ & $(+, 4, 15)$ & $132$\\
\hline
\hline
$327$ & $(+, 2, 17)$ & $145$\\
\hline
$331$ & $(+, 10, 11)$ & $145$\\
\hline
\hline
$389$ & $(-, 5, 13)$ & $177$\\
\hline
$399$ & $(+, 10, 13)$ & $177$\\
\hline
\hline
$395$ & $(-, 2, 17)$ & $179$\\
\hline
$403$ & $(+, 9, 14)$ & $179$\\
\hline
\hline
$409$ & $(-, 3, 16)$ & $186$\\
\hline
$417$ & $(+, 7, 16)$ & $186$\\
\hline
\hline
$421$ & $(-, 4, 15)$ & $192$\\
\hline
$427$ & $(+, 3, 19)$ & $192$\\
\hline
\hline
$431$ & $(-, 5, 14)$ & $197$\\
\hline
$439$ & $(+, 5, 18)$ & $197$\\
\hline
\hline
$449$ & $(-, 8, 11)$ & $206$\\
\hline
$451$ & $(-, 3, 17)$ & $206$\\
\hline
\hline
$453$ & $(+, 4, 19)$ & $204$\\
\hline
$455$ & $(-, 7, 12)$ & $204$\\
\hline
\hline
$461$ & $(-, 1, 20)$ & $210$\\
\hline
$463$ & $(+, 1, 21)$ & $210$\\
\hline
$469$ & $(+, 12, 13)$ & $210$\\
\hline
\end{tabular}
\end{center}
\end{minipage}

\bigskip

As stated in Theorem \ref{thm:Alex},
the Alexander polynomials of hyperbolic Berge knots of types VII and VIII
are not those of torus knots and satellite knots.
However we do not know the following:

\begin{qu}\label{qu:identify}
Do the Alexander polynomials of hyperbolic Berge knots of types VII and VIII
identify the knots ?
\end{qu}

From the table above, 
the degrees of the Alexander polynomials do not identify 
hyperbolic Berge knots of types VII and VIII completely.

\section{Final Remark}\label{sec:rem}
We give a remark about
a relation between our method and 
K.~Ichihara, T.~Saito and M.~Teragaito \cite{IST}
They show a formula of 
the Alexander polynomial of a doubly primitive knot in $S^3$.
Let $K$ be a doubly primitive knot in $S^3$.
Suppose that $(K; p)$ with $p\ge 2$ is a lens space $L(p, q)$.
Then the dual knot of $K$ in $L(p, q)$ can be expressed by
using a genus two Heegaard splitting of $L(p, q)$ and
a certain parameter $k$ ({\it Saito parameter}) with $1\le k\le p-1$.
We denote the knot by $K(L(p, q), k)$.

\medskip

For $i\in \mathbb{Z}$, let $\Psi(i)$ be the integral lift
of $i\overline{q}\ (\mathrm{mod}\ \! p)$
such that $1\le \Psi(i)\le p$, and
$$\Phi(i)=\sharp 
\{ j\ |\ \Psi(j)<\Psi(i)\quad \mbox{and}\quad 1\le j\le k-1\},$$
where $\sharp \{ \cdot \}$ implies the cardinality.

\begin{thm}\label{thm:IST}
{\rm (Ichihara, Saito and Teragaito \cite[Theorem 1.1]{IST})}
The Alexander polynomial of a doubly primitive knot $K$,
whose dual knot is $K(L(p, q), k)$, is 
$${\it \Delta}_K(t)\doteq
\frac{t-1}{t^k-1}\cdot \sum_{i=0}^{k-1}t^{\Phi(i)p-\Psi(i)k}.$$
\end{thm}

\begin{lm}\label{lm:IST}
{\rm (K.~Ichihara, T.~Saito and M.~Teragaito \cite[Lemma 2.1]{IST})}
$\gcd(k, p)=1.$
\end{lm}

\noindent
{\bf Proof}\ 
We give an alternative proof here.
By Theorem \ref{thm:IST}, we have
$$\left(\sum_{i=0}^{k-1}t^i\right)\cdot {\it \Delta}_K(t)\doteq
\sum_{i=0}^{k-1}t^{\Phi(i)p-\Psi(i)k}.$$
Suppose that $d=\gcd(k, p)\ge 2$.
Substituting $t=\zeta_d$ to the equation above, we have
$$0=k\ne 0.$$
This is a contradiction.
\qed

\bigskip

By the surgery formula II (Lemma \ref{lm:surgery2} (1)),
the Reidemeister torsion of $M=(K; p)$ is
\begin{eqnarray*}
\tau^{\psi_d}(M) & \doteq &
\frac{\zeta_d-1}{\zeta_d^k-1}\cdot 
\frac{\zeta_d^{k^2\overline{q}}-1}{\zeta_d^{k\overline{q}}-1}\cdot 
\frac{1}{(\zeta_d-1)^2}
\medskip \\
& \doteq &
\frac{\zeta_d^{k^2\overline{q}}-1}
{(\zeta_d-1)(\zeta_d^k-1)(\zeta_d^{k\overline{q}}-1)},
\end{eqnarray*}
where $d\ge 2$ is a divisor of $p$.

\begin{lm}\label{lm:IST2}
We have 
$$q\equiv \pm k^2\quad (\mathrm{mod}\ \! p)$$ 
and
$$\tau^{\psi_d}(M)\doteq 
(\zeta_d^k-1)^{-1}(\zeta_d^{\overline{k}}-1)^{-1}.$$
\end{lm}

\noindent
{\bf Proof}\ 
By the Franz lemma (Lemma \ref{lm:Franz}), we have
$$k^2\overline{q}\equiv
\pm 1\quad \mbox{or}\quad 
\pm k\quad \mbox{or}\quad 
\pm k\overline{q}\quad (\mathrm{mod}\ \! p).$$
From one of the latter two cases, we have
$$\tau^{\psi_d}(M)\doteq (\zeta_d-1)^{-2},$$
which is the Reidemeister torsion of $L(p, \pm 1)$.
By results due to P.~Ozsv\'ath and Z.~Szab\'o \cite{OS} 
and M.~Tange \cite{Ta}, $K$ is the unknot or the trefoil.
The knots satisfy the conditions.
From the first case, we have the result.
\qed

\bigskip

By Lemma \ref{lm:IST2},
$\{ \pm k, \pm \overline{k}\ (\mathrm{mod}\ \! p)\}$
has the same information as the Reidemeister torsions.
In \cite{ST}, 
$k\equiv -n\cdot \overline{m+n}\ (\mathrm{mod}\ \! p)$
has been obtained for the case $K=b^+(m, n)$.
The relations (\ref{eq:torb}) and (\ref{eq:torc}) 
can also be obtained from the facts.
M.~Tange gave a comment to the author that
since the Saito parameter $k$ is determined for every Berge knot,
to determine hyperbolicity of every Berge knot is not difficult
by the same way as the present paper.

\bigskip 

\noindent
{\bf Acknowledgement}\ 
The author would like to thank to
Kenneth Baker, Toshio Saito, and Masakazu Teragaito
who kindly answered elementary questions by the author,
and to Ruifeng Qiu and Motoo Tange
for giving him useful comments.
He also express his special gratitude
to Yuichi Yamada for consenting readily to divide 
the present work from our joint work \cite{KY}.

\bigskip


{\small
\par
Teruhisa KADOKAMI \par
Department of Mathematics, East China Normal University, \par
Dongchuan-lu 500, Shanghai, 200241, China \par
{\tt mshj@math.ecnu.edu.cn, kadokami2007@yahoo.co.jp} \par
}

\end{document}